\documentclass[12pt,reqno]{amsart}
\usepackage{color,latexsym}
\usepackage{listings}
\usepackage{graphicx}
\usepackage{subfig}

\newcommand{\e}{\epsilon}

\newcommand{\om}{\omega}
\newcommand{\Om}{\Omega}

\numberwithin{equation}{section}
\newcommand{\figref}[1]{Figure \ref{#1}}
\newtheorem{theorem}{Theorem}[section]
\newtheorem{defn}[theorem]{Definition}

\newtheorem{remark}[theorem]{Remark}
\newtheorem{prop}[theorem]{Proposition}

\begin{document}

\title[  ]
{Pseudorandom Numbers For Conformal Measures}

\author[  ]
{Manfred Denker$^1$, Jinqiao Duan$^2$ and Michael McCourt$^2$ \\ \\
1. Mathematics Department\\ Pennsylvania State University State\\
College, PA 16802, USA\\
  \emph{E-mail}:  denker@math.psu.edu   \\
\\
2. Department of Applied Mathematics\\ Illinois Institute of
Technology,  Chicago, IL 60616, USA \\
  \emph{E-mail}: duan@iit.edu    }

\date{April 17, 2009 (Revised version)}




\begin{abstract}
We propose a new algorithm for generating pseudorandom
(pseudo-generic) numbers of conformal measures of a continuous map
$T$ acting on a compact space $X$ and for a H\"older continuous
potential $\phi:X\to \mathbb R$. In particular, we show that this
algorithm  provides  good approximations to generic points for
hyperbolic rational functions   of degree two and the potential
$-h\log |T'|$, where $h$ denotes the Hausdorff dimension of the
Julia set of $T$.
\end{abstract}

\maketitle


\section{Introduction}\label{s1}

Conformal measures for rational maps were introduced by Sullivan (\cite{Su})
in 1983 following ideas of Patterson (\cite{Pa}) for the  case of limit sets of
Fuchsian groups. Existence and uniqueness of such measures has been shown in
\cite{DU} for a wide class of rational maps. These measures are in general
singular with respect to Lebesgue measure and have no explicitly computable
distribution function. There are a few papers dealing with the numerical
computation of these (mostly fractal) measures (e.g. \cite{Del}), but there is
no work done concerning the construction of generic points according to the
following definition.

\begin{defn} Let $(X,T)$ be a continuous dynamical system on a compact space
  $X$ and let $\nu$ be a $T$-invariant  probability measure on the Borel field of
  $X$. A point $x\in X$ is called generic if for every continuous function
  $h\in C(X)$
\begin{equation}\label{generisch}
 \lim_{n\to\infty} \frac 1n\sum_{k=0}^{n-1} h(T^k(x))= \int h d\nu.
\end{equation}
\end{defn}
The convergence rate  in this theorem can be arbitrarily slow, as
discussed in \cite{Krengel}.  It depends on the function $h$. A
probabilistic error bound can be obtained from the central limit
theorem or large deviation results; see \cite{Denker} for a
general description of the central limit theorem problem in
dynamical systems, and in particular \cite{Denker2}
 about this issue for rational functions.

 Accordingly, we call a point $y \in X$ pseudo-generic  for
 $\nu$ if the
 equation (\ref{generisch}) holds up to some prescribed  precision or
 error. The construction
 of pseudorandom numbers by the linear congruential method is also based on
 the iteration theory of maps of the interval. These points are as well
 pseudo-generic, hence one may use the notion of pseudorandom numbers as
 well in the present situation.

The aim of this
 note is to define and analyze an algorithm for computable pseudorandom
 points. In many applications, a sequence generated by the iteration of a
 pseudo-generic point will produce points which can be viewed as
 asymptotically independent realizations of independent identically
 distributed random variables. This follows whenever the map is a weakly
 dependent sequence of random variables.

Let $X$ be a compact metric space and $T:X\to X$ be continuous.
A conformal measure $m$ for a continuous function $\phi\in C(X)$ is a
probability measure satisfying
\begin{equation}\label{konform}
m(T(A))=\int_A \exp[\phi(x)] m(dx)
\end{equation}
for every measurable set $A$ with the property that $T$ restricted to $A$ is
invertible. This definition is equivalent to the requirement that
$$ \int g(x) m(dx)= \int_{A} g(T(x))\exp[\phi(x)] m(dx),$$
for every bounded continuous function $g$ and any measurable set $A$ such that
$T$ is invertible on $A$ and the support of $g$ is contained in $T(A)$. We
shall call (\ref{konform}) the conformal equation. Examples for such measures
are provided by rational maps (see \cite{DU} among others) or self-similar
measures on fractal sets (see \cite{Kigami} among others).

Given an invariant measure $\mu$, equivalent to a conformal measure $m$, one
can construct a pseudo-generic point by the method of least square estimation,
i.e. by minimizing
$$ \sum_{g\in \mathfrak G} \left(\frac 1n\sum_{k=0}^{n-1} g(T^k(x))-\int g
  d\mu\right)^2$$
over a suitable subclass $\mathfrak G$. However, the integral involved is not
  known and has to be computed by other means. It can be approximated in some  cases using the Perron-Frobenius
  operator
$$ Ph(x)=\sum_{T(y)=x} h(y)\exp[-\phi(y)]$$
and the projection to the eigenspace of the maximal eigenvalue of this
operator. Here we follow another approach using a discretized version of the conformal equation
together with a least square estimate. In this way no integral or calculation
of eigenspaces is involved in
the algorithm. The algorithm is explained and analyzed in section 2 in general
terms.

Sections 3 and 4  are devoted to the special case of hyperbolic rational functions
of degree two. We demonstrate how the algorithm is implemented
in this case. Other cases can be investigated in a similar way.

The algorithm requires the computation
of points  in $X$ and the density of the invariant measure $\mu$ with respect
to the conformal measure $m$ at specific points. Hyperbolic rational maps on the Riemann sphere $S^2=\overline{\mathbb
  C}$ are characterized by the requirement that their Julia sets
$J(T)$ does
not contain parabolic or critical points (see \cite{Beardon}). It is known that
these maps have a unique conformal measure $m$ for every H\"older-continuous
potential $\phi:J(T)\to \mathbb R$ (\cite{Su}), and that there is a unique equivalent,
ergodic and $T$-invariant probability measure $\mu$. This property will guarantee
the the pseudo-generic points are approximating integrals with respect to the
invariant measure and that the Perron-Frobenius operator can be used to find
the density $d\mu/dm$ at specific points. Unfortunately, this operator
requires to calculate the Hausdorff dimension of $J(T)$ which we do here
numerically.  Moreover, since repelling periodic
points are dense in $J(T)$ we are able to construct dense sets of points in
the Julia set.

It is important to remark that the algorithm needs precise numerical
calculations. We discuss this issue in Section 4.

\section{Least squares and the conformal equation}

In this section we describe an algorithm leading to explicitly computable
pseudorandom points for a dynamical systems. We start listing the assumptions
we impose to hold: Let $T:X\to X$ be a continuous map on some compact metric
space $X$ with metric $d(\cdot,\cdot)$ and let $\phi\in C(X)$ be a continuous
function. Assume the following conditions to hold:
\begin{itemize}
\item[A.]
\begin{enumerate}
\item There exists a unique conformal measure $m$ for $T$ and $\phi$.
\item There exists a unique ergodic, $T$-invariant measure $\mu$ which is
  equivalent to $m$.
\item The Radon-Nikodym derivative $f(x)=\frac{dm}{d\mu}$ has a continuous
  version defined on $X$ with modulus of continuity $\omega:\mathbb R\to
  \mathbb R$.
\end{enumerate}
\item[B.] For each $n,p\in \mathbb N$, $p\ge n$, there are finite sets
  $X_{n,p}\subset X$, finite sets ${\mathfrak A}_n$ of measurable subsets of
  $X$ and continuous functions $g_A\in C(X)$ ($A\in {\mathfrak A}_n$) such that
\begin{enumerate}
\item Every set $A\in {\mathfrak A}_n$ satisfies $A\in \sigma({\mathfrak
    A}_{n+1})$ (i.e. is a union of elements in ${\mathfrak A}_{n+1}$) and
    every point in $X$ lies in at most $a^*$ elements form ${\mathfrak A}_n$,
    where $a^*$ is independent of $n$.
\item $d_n:= \sup_{A\in {\mathfrak A}_n} \max\{ \mbox{diam}(A),
  \mbox{diam}(T(A))\} \to 0$ as $n\to\infty$.
\item $\mbox{supp}(g_A)\subset T(A)$ and $0\le g_A\le 2$.
\item The sigma fields $\sigma(\{ g_A: A\in {\mathfrak A}_n\})$ and $\sigma(\{
  g_A\circ T\cdot 1_A: A\in {\mathfrak A}_n\})$ generate the Borel field of $X$
  as $n\to\infty$. Each $g_A$ with $a\in {\mathfrak A}_n$ is approximated
  arbitrarily close by a linear combination of functions in ${\mathfrak A}_{n+l}$ for $l\ge
  1$ sufficiently large.
\item For each $n$, $X_{n,p}\subset X_{n,p+1}$ and
$$ D_{n,p}=\sup_{x\in X} \inf_{y\in X_{n,p}} \max_{1\le i\le p}
d(T^i(x),T^i(y))\to 0$$
as $p\to \infty$.
\end{enumerate}
\end{itemize}

The algorithm for fixed $n\in \mathbb N$ assumes the existence of $m$, $\mu$,
${\mathfrak A}_n$, $g_A$ ($A\in {\mathfrak A}_n$) and the sets $X_{n,p}$. It
proceeds as follows:
\begin{enumerate}
\item Choose $z_A\in A$ for $A\in {\mathfrak A}_n$ and compute $f(z_A)$ and $f(T(z_A))$.
\item Let $p=n$. For $x\in X_{n,p}$ compute $\beta_n^2(x)$ by
\begin{equation}\label{algorithm} \sum_{A\in {\mathfrak A}_n} \left( f(T(z_A))\sum_{k=0}^{p-1}
  g_A(T^k(x)) -f(z_A)\sum_{k=0}^{p-1}g_A(T^{k+1}(x))e^{\phi(T^k(x))}\right)^2.
\end{equation}
\item If $\min_{x\in X_{n,p}} \frac 1{p^2} \beta_n^2(x)\le 5a^*\omega(d_n)$
  stop and go to step 4, if
not set $p=p+1$ and continue with step 2.
\item Let $\beta_n^2=\min_{x\in X_{n,p}} \beta_n(x)$. Choose $x_n^*\in X_{n,p}$
minimizing this expression, i.e.
$$ \beta_n^2(x_n^*)= \beta_n^2.$$
\end{enumerate}

\begin{remark} (1) The algorithm requires to apply several subroutines
explained by examples in the following sections:

- Calculation of the sets $X_{n,p}$ and the distances $D_{n,p}$ for each fixed
  $n$.

- Calculation of the sets ${\mathfrak A}_n$, $d_n$ and $a^*$.

- Calculation of the functions $g_A$.

- Calculation of the density at points $z_A$ and $T(z_A)$.

\noindent (2) There is a simpler algorithm which may not work in general but
is easier to implement and used in later sections. In this case all sets
$A$ have no mass on their boundaries (because these will be a finite union of
points and the conformal measure has no atoms).

The simplification puts all functions $g_A$ to indicator functions, $g_A=1_{T(A)}$ and uses only one set
$X_{n,n}$.
\end{remark}

We need to show that the algorithm stops eventually and that the resulting
points $x_n^*$ are pseudorandom.
This will be accomplished in the following two propositions.

\begin{prop}
The algorithm stops eventually.
\end{prop}

\noindent{\it Proof.}  For $x\in X_{n,p}$ we have
\begin{eqnarray*}
&&\sqrt{\frac 1{p^2}\beta_n^2(x)} \\
&=&\left[\sum_{A\in {\mathfrak A}_n} \biggl(
f(T(z_A))\sum_{k=0}^{p-1}
 g_A(T^k(x))
 -f(z_A)\sum_{k=0}^{p-1}g_A(T^{k+1}(x))e^{\phi(T^k(x))}\biggr)^2\right]^{1/2}\\
&=&\left[\sum_{A\in {\mathfrak A}_n} \biggl(
f(T(z_A))\sum_{k=0}^{p-1}
 g_A(T^k(x)) - \int g_A dm\right. \\
&& + \left. \int_Ag_A(T(u)) e^{\phi(u)} m(du)
 -f(z_A)\sum_{k=0}^{p-1}g_A(T^{k+1}(x))e^{\phi(T^k(x))}\biggr)^2\right]^{1/2}\\
&\le&\left[\sum_{A\in {\mathfrak A}_n} \biggl(
f(T(z_A))\sum_{k=0}^{p-1}
 g_A(T^k(x)) - f(T(z_A)) \int g_A d\mu\right. \\
&& + \left. f(z_A)\int_Ag_A(T(u)) e^{\phi(u)} \mu(du)
 -f(z_A)\sum_{k=0}^{p-1}g_A(T^{k+1}(x))e^{\phi(T^k(x))}\biggr)^2\right]^{1/2}\\
&+&
  \left[\sum_{A\in {\mathfrak A}_n} \biggl(f(T(z_A)) \int g_A d\mu - \int g_A fd\mu\biggr)^2\right]^{1/2}\\
&+& \left[\sum_{A\in{\mathfrak A}_n}\biggl(f(z_A)\int_Ag(T(u))
e^{\phi(u)} \mu(du)-\int_Ag_A(T(u)) e^{\phi(u)}
 f(u) \mu(du)\biggr)^2\right]^{1/2}\\
&\le&\biggl[\sum_{A\in {\mathfrak A}_n} \biggl(
f(T(z_A))\sum_{k=0}^{p-1}
 g_A(T^k(x)) - f(T(z_A)) \int g_A d\mu \\
&& + f(z_A)\int_Ag_A(T(u)) e^{\phi(u)} \mu(du)
 -f(z_A)\sum_{k=0}^{p-1}g_A(T^{k+1}(x))e^{\phi(T^k(x))}\biggr)^2\biggr]^{1/2}\\&+& 4a^*\omega(d_n).
\end{eqnarray*}
Since $\mu$ is ergodic there exists a generic point for $\mu$. Let
$z$ denote such a point. Choose $x\in X_{n,p}$ such that
$d(T^i(z),T^i(x))\le D_{n,p}$, for each $1\le i\le p$. Then for
any $g_A$ by the triangle inequality
\begin{eqnarray*}&&|\int g_A d\mu - \frac 1p\sum_{k=0}^{p-1} g_A(T^k(x))|\\
&\le & |\int g_A d\mu - \frac 1p\sum_{k=0}^{p-1} g_A(T^k(z))| + \omega_{g_A}(D_{n,p})
\end{eqnarray*}
which converges to 0 as $p\to \infty$, where $\omega_h$ denotes the modulus of
continuity for the function $h$.
Likewise
\begin{eqnarray*}&&
|\int_Ag_A(T(u)) e^{\phi(u)} \mu(du)- \frac 1p\sum_{k=0}^{p-1}
 e^{\phi(T^k(x))}g_A(T^{k+1}(x))1_A(T^k(x))|\\
&\le & |\int_Ag_A(T(u)) e^{\phi(u)} \mu(du)- \frac 1p\sum_{k=0}^{p-1}
 e^{\phi(T^k(x))}g_A(T^{k+1}(x))1_A(T^k(z))|+\omega_{e^{\phi}g_A\circ T},\\
\end{eqnarray*}
which tends to 0 as well as $p\to \infty$ (note that the support
of $g_A\circ T$ is inside of $A$, so that the discontinuity of $A$
is of no relevance). The indicator of A has to be kept here. If
$g_A$ is the indicator of $T(A)$, then $g_A\circ T(x)$ is one if
and only if $x$ is a preimage of a point in $T(A)$. But there are
more than the points in $A$ as preimages and we want only those in
$A$.

Since ${\mathfrak A}_n$ is a finite set, not changing with $p$, we
see that
$$ \limsup_{p\to\infty} \frac 1{p^2}\beta_n^2\le 4a^*\omega(d_n),$$
which implies that the algorithm stops eventually.\qed

\begin{prop}
Let $x_n^*$ and $p=p(n)\ge n$ ($n\in\mathbb N$) be constructed
according to the algorithm, and assume  that A and B hold. Then
for every continuous function $g$ we have
$$ \lim_{n\to\infty} \frac 1{p(n)} \sum_{k=0}^{p(n)-1} g(T^k(x_n^*))=\int g
d\mu.$$\end{prop}

\noindent{\it Proof.} Define
$$ \nu_n= \frac 1{p(n)} \sum_{k=0}^{p(n)-1} \delta_{T^k(x_n^*)},$$
where $\delta_z$ denotes the point mass in $z\in X$. Then $\{\nu_n:n\in
\mathbb N\}$ is relatively compact in the weak topology of measures. Let $\nu$
be an accumulation point.
Then $\nu$ is an invariant measure and for $A\in {\mathfrak A}_n$
\begin{eqnarray*}
&&| \int g_A fd\nu_{n+l}- \int_Ae^{\phi(u)} f g_A(T(u)) \nu_{n+l}(du)|\\
&\le&
   \sum_{g_B:b\in {\mathfrak B}\subset {\mathfrak A}_{n+l}}| f(T(z_B)) \int g_B
   d\nu_{n+l}- f(z_B)\int_Be^{\phi(u)} g_B(T(u)) \nu_{n+l}(du)| +o_l(1)\\
&\le& \beta_{n+l} +o_l(1).
\end{eqnarray*}
Therefore, letting $l \to\infty$ along a suitable subsequence,  $d\widetilde m= f
d\nu$ satisfies the conformal equation, hence is conformal.  Since $m$ is
unique as a conformal measure it equals $\widetilde m$. Moreover, $\nu$ is
equivalent to $m$, and must be equal to $\mu$, since the latter is unique as
well. This shows that $\nu_n$ converges weakly to $\mu$, proving that $x_n^*$
is a sequence of pseudorandom points. \qed

\section{Conformal measure on a Julia set}  \label{s2}

We describe a specific example in this section, which will be used to show how
the general algorithm can be applied.

Consider the  rational map $T:C \to C$   defined by
$T(z)=z^2+\frac18$. Its Julia set $J(T)$ is a bounded compact set
in $C$ (under the induced (Euclidean) topology).

A Julia set $J(T)$ is the closure of the set of repelling periodic
points and inverse images of $T^n$, ($n\ge 1$) are dense as well
(\cite{Beardon}). Thus the repelling periodic points and its preimages are dense inside the
Julia set $J(T)$, and every point in   $J(T)$ can be realized as a
limit of some sequence of preimages of each repelling periodic point.

   The set of
  points $y$ such that the iterates of $y$ under $T$, $T^2$, $T^3$ and so on
  eventually hit a fixed repelling periodic point $z$ is dense in the Julia set. Therefore, it is possible to construct many
  points in the Julia set by taking preimages. This can be done in different
  ways. The most convenient is to calculate inverse branches $f_1$ and $f_2$
  (for a quadratic polynomial) as maps defined
  on the Julia set.  Then we can iterate all possible finite combinations $f_1\circ
  f_2\circ f_1\circ f_1...$, where the sequence of $f_1$ and $f_2$'s are
  arbitrary choices. This is the naive way, since the computation creates too
  many data. We need for later purpose a certain depth of the iteration
  procedure, much longer than the forward iteration done later. In order to
  accomplish this one takes random choices of the two maps over a long string
  of iterations. This gives one point in the Julia set and one needs to
estimate the errors in this calculation.

We now discuss the discretization, i.e., a random mesh on the
Julia set.

Take a repelling periodic point $z_0$, say a fixed point:
$T(z_0) =z_0$. Let us
discretize the Julia set $J(T)$, i.e., generate a random mesh or
random lattice $S$ over it, as its computational representation.

The inverse of $T$ has two analytical branches, $f_1$ and $f_2$.
We backward iterate $T$ but select the inverse branches randomly.
Let $l $ be a large positive integer. Define a sample space
$$
\Omega = \{1, 2  \}^l = \{ \omega=(\om_1, \om_2, \cdots, \om_l)\;
: \; \om_i = 1 \; \mbox{or}\; 2 \}.
$$
Starting from $z_0$, a random backward iteration of $l$ steps of
$T$ can be represented as
$$
f_{\om_l}  \cdots f_{\om_2} f_{\om_1} (z_0),
$$
where $\om_1,...,\om_l$ are chosen randomly with equal probability.
Let $\Om_0$ be a randomly chosen set of $\om_1,...,\om_l$. Its cardinality is
the the size of the random
mesh.
Now we define a random mesh of the Julia set $J(T)$ defined by $\Om_0$ as
\begin{eqnarray}\label{lattice}
S = \{z:\; z=f_{\om_l}  \cdots f_{\om_2} f_{\om_1} (z_0),
\om=(\om_1, \om_2, \cdots, \om_l) \in \Om_0 \}.
\end{eqnarray}

We take the Borel sets $A$ to be small balls in $J(T)$ centered around
some points in $S$. Let us take a representative subset $S_0$ of
$S$ and take $A$ as balls centered on points in $S_0$.   This is a
finite family  of balls and it is arranged to cover $J(T)$:
\begin{eqnarray}  \label{family}
 \mathcal{A} :=\{A=B(z^*, \delta) : \; \mbox{ball  with center} \;  z^*\in S_0  \;
\mbox{and radius}\;  \delta>0     \}.
\end{eqnarray}

Let $h$ denote the Hausdorff dimension of the Julia set $J(T)$.
We shall consider the conformal measure $m$ associated to the potential $h\log
|T'|$, which is a well defined Lipschitz continuous function on $J(T)$, since
$T$ is hyperbolic.
The transfer operator (Perron-Frobenius operator) for the the map $T$ and the
potential is
$P:  C(J(T))  \to  C(J(T))$   defined as
\begin{eqnarray}
P g(z) = \sum_{y \in T^{-1}(z)} g(y)|T'(y)|^{-h}.
\end{eqnarray}
Iteration yields
\begin{eqnarray}
P^n g(z) = \sum_{y \in T^{-n}(z)} g(y)|(T^n)'(y)|^{-h}.
\end{eqnarray}
It is known that
 there exists a unique
conformal measure with respect to this potential, always denoting it $m$. Moreover, $T$ is
ergodic and has a unique finite invariant measure  $\mu$ (on
$J(T)$)  that is equivalent to the conformal measure $m$.

\medskip
  The Hausdorff dimension of the Julia set can
 be calculated as follows. It is known that
$$ \sum_{n=0}^\infty \sum_{T^n(y)=x} |(T^n)'(y)|^s$$
converges for $s< -h$ and diverges for $s>-h$. We shall use a slight variant of
this fact to determine $h$:
$$\sum_{T^n(y)=x} |(T^n)'(y)|^s$$
converges only for $s=-h$.
The Hausdorff dimension $h$ is approximately $=1.00735$ for the map $z\mapsto
z^2+\frac 18$.

\medskip

The starting point for the optimization is the defining equation for a conformal
measure \eqref{konform}. We evaluate this equation for  balls  $A\subset J(T)$.
It is known that $m$ has no atoms for our special quadratic map considered here,
since the Julia set is a Jordan curve. Thus $m(\partial B)=0$ for open sets   and we can construct
pseudo-generic points from the equation
\begin{eqnarray}
m(T(A)) = \int_A |T'|^h dm
\end{eqnarray}
directly, where $A$ are balls.

Since $\mu \sim m$, by the Radon-Nikodym theorem
$\frac{d\mu}{dm} = f(z)$ (density of $\mu$ w.r.t. $m$), $\mu
(TA)=\int_{TA} f dm $ or
$$
m (TA)= \frac1{f(Tz^{**})} \mu(TA)=\int_A |T'|^h dm = \frac1{f(z^*)} \int_A |T'|^h d\mu
$$
  for some point $z^*$ and $z^{**}$ by the intermediate value theorem since the density
  is continuous. These equations hold approximately for all $z$
  replacing $z^*$ and $z^{**}$ if the
sets $A$ and $T(A)$ are small enough.

Thus, on small balls $A$ in the Julia set, we have
\begin{eqnarray}
m(T(A)) & =&    \frac1{f(T(z^*))} \int_{T(A)} d\mu   \nonumber\\
& =& \frac1{f(T(z^*))} \lim_{n\to \infty} \frac1{n}
\sum_{k=0}^{n-1} 1_{T(A)} (T^k(z)), \\
\int_A |T'|^h dm &=& \frac1{f(z^*)} \int_A |T'|^h d\mu  \nonumber
\\
 &=&\frac1{f(z^*)}\lim_{n\to \infty} \frac1{n} \sum_{k=0}^{n-1}
1_A (T^k(z)) \; |T'(T^k(z))|^h,
\end{eqnarray}
where $z\in J(T)$ is in a full $\mu-$measure subset in $J(T)$.

 Therefore, we have the fundamental equation
\begin{eqnarray}  \label{eqn}
 \frac1{f(T(z^*))}  \lim_{n\to \infty} \frac1{n}
\sum_{k=0}^{n-1} 1_{T(A)} (T^k(z)) =  \frac1{f(z^*)}  \lim_{n\to
\infty} \frac1{n} \sum_{k=0}^{n-1} 1_A (T^k(z)) \; |T'(T^k(z))|^h,
\end{eqnarray}

We will find such a $z$ approximately by the method of least squares, i.e.
finding the minimizer for
\begin{eqnarray}  \label{optimal}
\min_{z \in S}  \sum_{A \in \mathcal{A}} \left |\frac1{f(T(z^*))}
\frac1{n} \sum_{k=0}^{n-1} 1_{T(A)} (T^k(z)) -\frac1{f(z^*)}
\frac1{n} \sum_{k=0}^{n-1} 1_A (T^k(z)) \; |T'(T^k(z))|^h \right|^2,
\end{eqnarray}
where -more precisely- $z^*$ depends on the corresponding ball $A$.

For calculating  the density $f(z^*)$ we use the
  transfer  operator and the  well know equation
\begin{eqnarray}\label{PerFrob}
f(z)=\lim_{n\to\infty}P^n 1(z)=\lim_{n\to \infty} \sum_{y \in T^{-n}(z)}  |(T^n)'(y)|^{-h}
\end{eqnarray}
We choose points $z^*\in S_0$ in the discretization step and open balls around
these points as choices of the sets $A$. The equation (\ref{PerFrob}) is used
to
calculate the densities at $z^*$ and $T(z^*)$.

\bigskip

We check numerically whether the minimizer in
the optimization problem is pseudo-generic.

For any continuous and bounded function on $J(T)$, and for any
generic point $z$ in $J(T)$, we should have
\begin{eqnarray}  \label{ergodic}
\int g d \mu = \lim_{n\to \infty}  \frac1{n}  \sum_{k=0}^{n-1}
g(T^k z), \; a.e.-\mu,
\end{eqnarray}
for $g \in L^1(J(T))$. In fact, in the next section, we test this
for the function $g=|z|$.
 For $z$ obtained in our numerical procedure, we
compute  $\lim_{n\to \infty} \frac1{n} \sum_{k=0}^{n-1} g(T^k
(z))$ for some large $n$. Repeated calculation for different sets
$S$,
 $S_0$, random choices  and $n$ will show that the average does not vary
 considerably. Choosing the backwards iteration  randomly and $S_0$ randomly
 may be seen as the analog of a seed in the construction of random numbers.

\section{Numerical results} \label{mike}

\subsection{Determining the Hausdorff Dimension $h$ for the Julia
Set}

\begin{figure}[h]
\includegraphics[width=7cm]{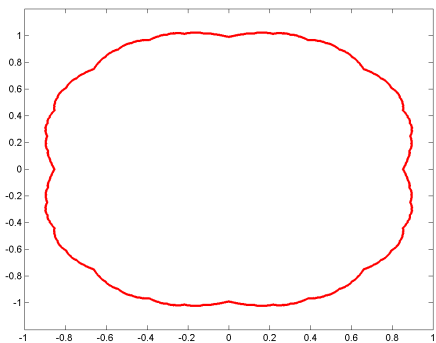}
\caption{\footnotesize{Julia Set of the mapping $T(z)=z^2+1/8$}}
\label{Julpic}
\end{figure}

Before we can hope to evaluate the Perron-Frobenius operator for a
generic point in $J(T)$, we must determine the appropriate
Hausdorff dimension.  There should be only one value $h$ which
allows for convergence of the limit described by \eqref{PerFrob}
to a value $f(z)\in(0,\infty)$.  This $h$ will be different for
each rational map but since we only consider one such map here,
$T(z)=z^2+\frac{1}{8}$, we need only determine the dimension of
the resulting fractal (Fig. \ref{Julpic}).

We begin by taking the only repelling fixed point of $T$ as a test value to determine $h$.
Below is a table which shows the sequence of $f_n(z)$ values (see \S\ref{fofz}) as $n\to\infty$.

\begin{figure}[h]
\begin{tabular}{|c|cc|cc|cc|r|}
\multicolumn{7}{c}{\textbf{Convergence of the Density Function}} \\ \hline
$z_0=.8356$&\multicolumn{2}{c}{$h=1.00$}&\multicolumn{2}{c}{$h=1.00735$}&\multicolumn{2}{c} {$h=1.01$} \\
&$f_n(z_0)$&$f_n/f_{n-1}$&$f_n(z_0)$&$f_n/f_{n-1}$&$f_n(z_0)$&$f_n/f_{n-1}$ \\ \hline
$n=3$&1.3029&1.0983&1.2922&1.0935&1.2884&1.0918 \\
$n=5$&1.4132&1.0299&1.3884&1.0250&1.3796&1.0232 \\
$n=10$&1.4865&1.0059&1.4245&1.0009&1.4028&0.9991 \\
$n=15$&1.5256&1.0051&1.4258&1.0000&1.3914&0.9982 \\
$n=20$&1.5644&1.0050&1.4258&1.0000&1.3789&0.9982 \\
\hline
\end{tabular}
\caption{\footnotesize{Table to experimentally determine the Hausdorff dimension.}}
\end{figure}

It is clear that the approximate value $h=1.00735$ allows for
convergence of the Perron-Frobenius operator, and that values too
small or too big yield unbounded or zero answers respectively. For
the remainder of this project we will approximate $h$ as such.
More rigorous discussion of Hausdorff dimension computation can be
found in \cite{Kigami}.
If $n$ increases the results  become better in generally. Stratistical
rigourous methods use the  Grassberger and Procaccia correlation dimension and
has been developed by Cutler and
  Dawson or Denker and Keller. The above approach is sufficient for hyperbolic maps.

\subsection{Creating the Computational Lattice}\label{algo}


Now that we have determined the appropriate Hausdorff dimension,
we  compute the $z^*$ which define the $A$ (the covering of the Julia set) and the lattice $S$.  Given an $m$, the $z^*$ are all the points for which $z_0=T^m(z^*)$, which can be computed directly without much difficulty.  The second block of code in Appendix A.1 generates the points $z^*$ using $z_0=.8536$.  Logically, there are $2^m$ points which generate the covering of the Julia Set, since there are $2^m$ pre-images in $T^{-m}(z_0)$.  An example of a covering of the Julia set is below.

\begin{figure}[h]
\includegraphics[scale=.6]{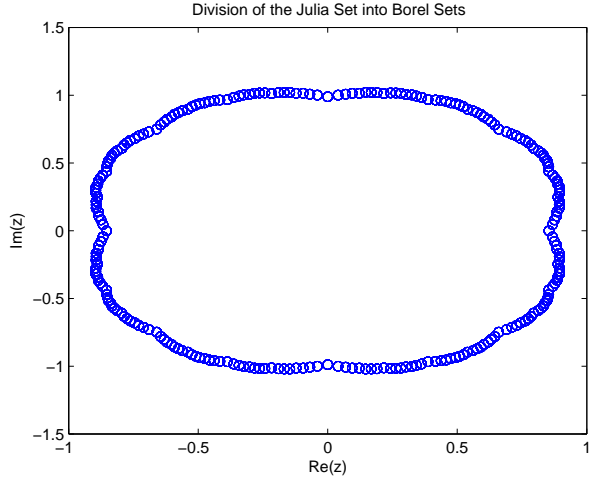}
\label{borelcover}
\caption{\footnotesize{For $m=8$ there are 256 $z^*$ points which center the circles $A\in\mathcal{A}$ which cover $J(T)$.}}
\label{graphbranch2}
\end{figure}

Each of the points in $S$ is determined by an inverse iteration backwards from a point randomly chosen from the $z^*$.  $|S|=\ell$, and each point is inverse iterated backwards $N$ times from the initial randomly chosen point.  At each inverse iteration only one pre-image is chosen and stored, since the other pre-image is of no consequence for determining the final lattice.  The appropriate Matlab code to create the lattice is the function \texttt{makelattice} which is found in Appendix A.2.

\subsection{Numerical Error Analysis in  Random Mesh Generation }

In \S \ref{algo}, we use the random mesh $S$ to discretize the
Julia set $J(T)$. So at least, we would need to make sure that
points in $S$ are (approximately) inside $J(T)$, when $\ell$ is large
enough. Numerical error here comes mainly from computing the
inverse branches $f_1=F_1, f_2=F_2$ of $T=z^2+\frac18$. Let us show that
this numerical process is stable, i.e., the total error is bounded
(\cite{Kincaid}).


Let $F$ denote either one of $F_1$ and $F_2$. The error analysis
for numerically iterating each function is similar. Let the
computer's unit roundoff error be, for example with double
precision, $ 10^{-16}$.  Let us ignore the error in computing the
initial repelling fixed point $z_0$.

Denote $ \tilde{F} (z_0)$ be the computed value of $F(z_0)$. Then
$$ | \tilde{F} (z_0) - F(z_0)| \leq |F(z_0)| \e_0, \; \e_0 \leq 10^{-16}. $$
Denote $z_1= F(z_0),  \tilde{z}_1= \tilde{F}(z_0)$. In the
following $\e_i$'s denote relative roundoff errors at various
steps of computation. Then
\begin{eqnarray*}
 | \tilde{F} (\tilde{z}_1) - F(z_1)|
 &=&  | \tilde{F} (\tilde{z}_1) - F(\tilde{z}_1) + F(\tilde{z}_1)- F(z_1)|   \\
 &\leq & | \tilde{F} (\tilde{z}_1) - F(\tilde{z}_1)| + |F(\tilde{z}_1)- F(z_1)|   \\
 &\leq &  |F(z_1)| \e_1+ |F'(z_1) (\tilde{z}_1-z_1)| \\
 &\leq & |F(z_1)| \e_1+ |F'(z_1)|\; |  F(z_0)| \e_0, \;  \; \e_1 \leq 10^{-16}.
\end{eqnarray*}
As long as $F$ and its  derivative is bounded in the (bounded)
Julia set, the right hand side of the above error estimate is
upper bounded by a constant multiplying $\e \leq 10^{-16}$. Thus
when we iterate $n$ times of $F$, the relative roundoff error
(\cite{Kincaid}) is approximately $n\e$. The random selection
between $F_1$ and $F_2$ does not change this error order.



\subsection{Efficiently Computing the Transfer Operator}\label{fofz}

Now that we have determined our computational lattice we must compute the
transfer operator for obtaining $f(z)$ which is the density as defined by the Radon-Nikodym derivative.
This is a computationally sensitive segment of the procedure because it requires
approximating a limit which becomes exponentially more expensive to compute.  We can make
this somewhat easier by noting that we consider only one mapping, and that $T'(z)$ is a linear function.
\begin{equation}\label{twoz}
 T(z)=z^2+\frac{1}{8} \quad\Rightarrow\quad T'(z)=2z
\end{equation}
The chain rule allow us to say the following,
\begin{equation}\label{chainrule}
 (T^n)'(z)=(T\circ T^{n-1})'(z)=T'(T^{n-1}(z))T'(T^{n-2}(z))...T'(T(z))T'(z).
\end{equation}

Below we substitute \eqref{twoz} and \eqref{chainrule} into \eqref{PerFrob} and use the chain rule
technique mentioned to simplify the evaluation of the transfer operator.  The final simplification occurs because
\begin{equation*}
|uv|=|u||v| \qquad \forall u,v\in\mathbb{C}
\end{equation*}
Note that the set $T^{-n}(z)$ includes all $2^n$ (possible
non-unique) values on the Julia set for which $T^n(T^{-n}(z))=z$.

\begin{align}
f(z)&=\lim_{n\to \infty} \sum_{y \in T^{-n}(z)}  |(T^n)'(y)|^{-h} \nonumber \\
&=\lim_{n\to\infty} \sum_{y \in T^{-n}(z)} |T'(T^{n-1}(y))...T'(y)|^{-h} \nonumber \\
&=\lim_{n\to\infty} \sum_{y \in T^{-n}(z)} |2T^{n-1}(y)...2y|^{-h} \nonumber \\
&=\lim_{n\to\infty} 2^{-hn} \sum_{y \in T^{-n}(z)} \left|\prod_{k=0}^{n-1}T^k(y)\right|^{-h} \nonumber \\
&=\lim_{n\to\infty} 2^{-hn} \sum_{y\in T^{-n}(z)} \left(\prod_{k=0}^{n-1}\left|T^k(y)\right|\right)^{-h}
\label{middletrans}
\end{align}

We can further simplify this by noting that we are not interested in the actual values
of $T^{k-n}(z)$ but rather only their moduli.  For each value $y\in T^{-n}(z)$ there is
a value $-y\in T^{-n}(z)$ since both the positive and negative square roots must be considered.
Because of this we need only calculate the contribution of half the pre-images to the summation
and then double it since $|y|=|-y|$.

If we let $T^{-n}_+$ denote only the positive square root pre-images of all $2^{n-1}$ pre-images
in $T^{-(n-1)}$ (note that $|T^{-(n-1)}|=|T^{-n}_+|$), we can modify \eqref{middletrans}
\begin{align}
f(z)&=\lim_{n\to\infty} 2^{-hn} \sum_{y\in T^{-n}(z)} \left(\prod_{k=0}^{n-1}\left|T^k(y)\right|\right)^{-h} \nonumber \\
&=\lim_{n\to\infty} 2^{-hn} \left[2\sum_{y\in T^{-n}_+(z)} \left(\prod_{k=0}^{n-1}\left|T^k(y)\right|\right)^{-h}\right] \nonumber \\
&=\lim_{n\to\infty} 2^{-hn+1}\sum_{y\in T^{-n}_+(z)} \left(\prod_{k=0}^{n-1}\left|T^k(y)\right|\right)^{-h}
\label{finaltrans}
\end{align}

\begin{figure}[h]
\includegraphics[width=10cm]{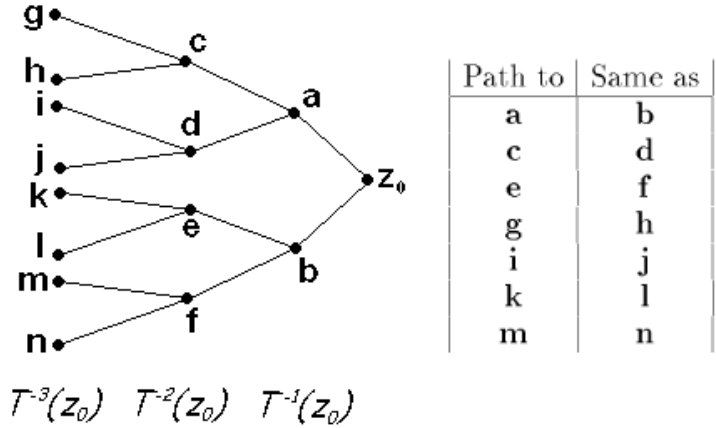}
\caption{\footnotesize{Mutliple pre-image paths yield the same product in
computing the density in $(\ref{finaltrans})$} }
\label{graphbranch1}
\end{figure}

In the Matlab code to execute this we take advantage of the
 fact that we need only calculate and store half the moduli
 needed by the same logic used before.  This is seen in graph theory
 and Figure \ref{graphbranch1} shows how we utilize the fact that several
 branches in the tree have the same product.

\medskip

Using this logic we wrote the function \texttt{densop} to approximate the
 limit in \eqref{finaltrans}  to $10^{-4}$ accuracy.  If the appropriate $h$
 value is not used, the program will likely fail to converge and run in perpetuity.  See Appendix A.3 for the Matlab implementation.

  For those who are interested in the speed of this algorithm,
ours is certainly not the fastest possible implementation.
Each time \texttt{densop} is called it recalculates pre-images
which may have already been determined.  In addition, new memory
is allocated in each iteration above as well as at the start of
each call to \texttt{densop}.  All the code in this project
is written to test the algorithms described above and to emphasize readability; this has resulted in a decrease in efficiency which will be the topic of a future project.

\subsection{Solving the Optimization Problem}

Recall the optimization problem we need to solve:
given $z^*$ and $\delta_A$ (to define the Borel Set $A=B(z^*,{\delta_A})$) and $n$
\begin{equation*}
\min_{z \in S}  \sum_{A \in \mathcal{A}}   \left| \frac1{f(Tz^*)}
\frac1{n} \sum_{k=0}^{n-1} 1_{T(A)} (T^kz) -\frac1{f(z^*)}
\frac1{n} \sum_{k=0}^{n-1} 1_A (T^kz)
\left|T'(T^kz)\right|^h\right|^2.
\end{equation*}
One thing to note is that the optimization is to occur on a discrete
lattice, $S$, thus we need only test a finite number of points, $\ell$,
to find the solution.  Another point of interest is that $T^k(z)$ is
 evaluated during construction of the lattice so no new function evaluations
 take place.  Also there is no need to set $n>\ell$ since $T^\ell(z)=z_0$
 for all $z\in S$.

To determine if a point $z$ is in $B(z^*,{\delta_A})$ we simply
test $|z-z^*|<\delta_A$, thus
\begin{equation}\label{indicator}
 1_A(z)=\left\{\begin{array}{cc}1&\quad|z-z^*|<\delta_A\\
                                0&\quad\mbox{else}\end{array}\right.
\end{equation}

Testing if $z$ is in $B(z^*,{\delta_{T(A)}})$ is more difficult.
To do so we state that
\begin{align}
 z&\in T(A), \nonumber \\
 T^{-1}(z)&\in A,
\label{preimage}
\end{align}
and therefore test whether either pre-image of $z$ is in $A$.  The result is that
\begin{equation}\label{pindicator}
 1_{T(A)}(z)=\left\{\begin{array}{cc}1&\quad|T^{-1}_\pm(z)-z^*|<\delta_A\\
                                    0&\quad\mbox{else}\end{array}\right.
\end{equation}
where $T^{-1}_\pm(z)$ is either the positive or negative preimage of $z$.


Appendix A.4 is the Matlab implementation of the optimization procedure.  This code simply runs through every point in $z\in S$ and returns the point which minimizes \eqref{optimal} for summations of a given length $n$.  It also returns the $\beta_n$ value described by \eqref{algorithm}.

\subsection{Testing the Pseudorandom Points}
In order to determine if the $z$ which satisfies \eqref{optimal} is a pseudorandom point we test its time average in equation \eqref{ergodic}.  We use the simple test function $g(x)=|x|$ for which the integral can be approximated deterministically; the lhs of \eqref{ergodic} is the average distance of points in $J(T)$ from the origin.  When all the pre-images for various values of $m$ are averaged together, we see below that the result approaches the limit $\int|z|d\mu\approx1.001379$.

\begin{figure}[h]\label{ensemavg}
\begin{tabular}{|c|c|}
\multicolumn{2}{c}{\textbf{Integral Limit}} \\ \hline
$m$&$\int|z|d\mu$ \\ \hline
1 & 0.853553 \\
5 & 0.990741 \\
10 & 1.001044 \\
15 & 1.001369 \\
20 & 1.001379 \\
21 & 1.001379 \\ \hline
\end{tabular}
\caption{\footnotesize{Table describing the approximate solution to \eqref{ergodic}.}}
\end{figure}



There are several factors which contribute to the quality of the solution and the complexity of the algorithm.  We will assume here that we already know $h$ and that $m=8$ is fixed which means that the Julia set is covered by $2^8=256$ balls and that $\delta_A=2^{1-m}$ is also fixed.  Assuming that the $z^*$ are already available (which is reasonable because the Borel sets are defined by $m$), the only parameters which affect the accuracy of the integral are
\begin{itemize}
\item $\ell$ - The number of elements in $S$, the computational lattice
\item $N$ - The depth of the inverse iteration on randomly chosen points $z^*$ to generate $S$
\begin{itemize}
\item Recall that for $z\in S$, $T^N(z)=z^*$ for exactly one of the $2^m$ $z^*$.
\end{itemize}
\item $n$ - The point at which we truncate the least squares limit
\item $\alpha$ - The number of pseudorandom trajectories used to evaluate the integral
\end{itemize}

We can now take a look at how changing these parameters individually affects the speed and accuracy of the algorithm.  For all these experiments we have generated 10 computational lattices (ie $\alpha=10$); for each lattice 1 pseudorandom point minimizes the least squares equation and that is the point whose trajectory we use to compute the integral.  $\mu$ is the experimental mean and $\sigma$ is the experimental standard deviation

\begin{figure}[h]
\centering
\begin{tabular}{cc}
\subfloat[Fixed $N=16000$, $n=100$.]{\label{fig:diffell}
\begin{tabular}{|c|ccc|}\hline
$\ell$&$\mu$&$\sigma$ ($\times 10^{-3}$)&Time \\
\hline
25&1.00136&0.5848&6 \\
50&1.00179&0.8568&11 \\
100&1.00149&0.8295&20 \\
200&1.00145&0.9279&42 \\
400&1.00122&0.9090&91 \\
800&1.00150&0.9652&205 \\
1600&1.00150&0.5960&413 \\
3200&1.00142&0.4890&826 \\ \hline
\end{tabular}}
&
\subfloat[Fixed $N=16000$, $\ell=100$.]{\label{fig:diffn}
\begin{tabular}{|c|ccc|}\hline
$n$&$\mu$&$\sigma$ ($\times 10^{-3}$)&Time  \\
\hline
25&1.00117&0.8050&13 \\
50&1.00160&0.6210&15 \\
100&1.00149&0.8295&20 \\
200&1.00116&0.6325&32 \\
400&1.00148&0.8654&70 \\
800&1.00147&0.5014&146 \\
1600&1.00136&0.5816&281 \\
3200&1.00145&0.5851&576 \\ \hline
\end{tabular}}
\end{tabular}

\subfloat[Fixed $n=100$, $\ell=100$.]{\label{fig:diffN}
\begin{tabular}{|c|ccc|}\hline
$N$&$\mu$&$\sigma$ ($\times 10^{-3}$)&Time  \\
\hline
1000&1.00109&3.264&12 \\
2000&1.00116&2.635&13 \\
4000&1.00153&1.732&14 \\
8000&1.00108&0.8890&16 \\
16000&1.00149&0.8295&20 \\
32000&1.00138&0.4505&30 \\
64000&1.00119&0.3765&48 \\
128000&1.00145&0.1743&85 \\ \hline
\end{tabular}}
\caption{The effect of varying $N$, $n$ and $\ell$.\label{fig:diffparam}}
\end{figure}

\begin{figure}[h]
\centering
\includegraphics[scale=.7]{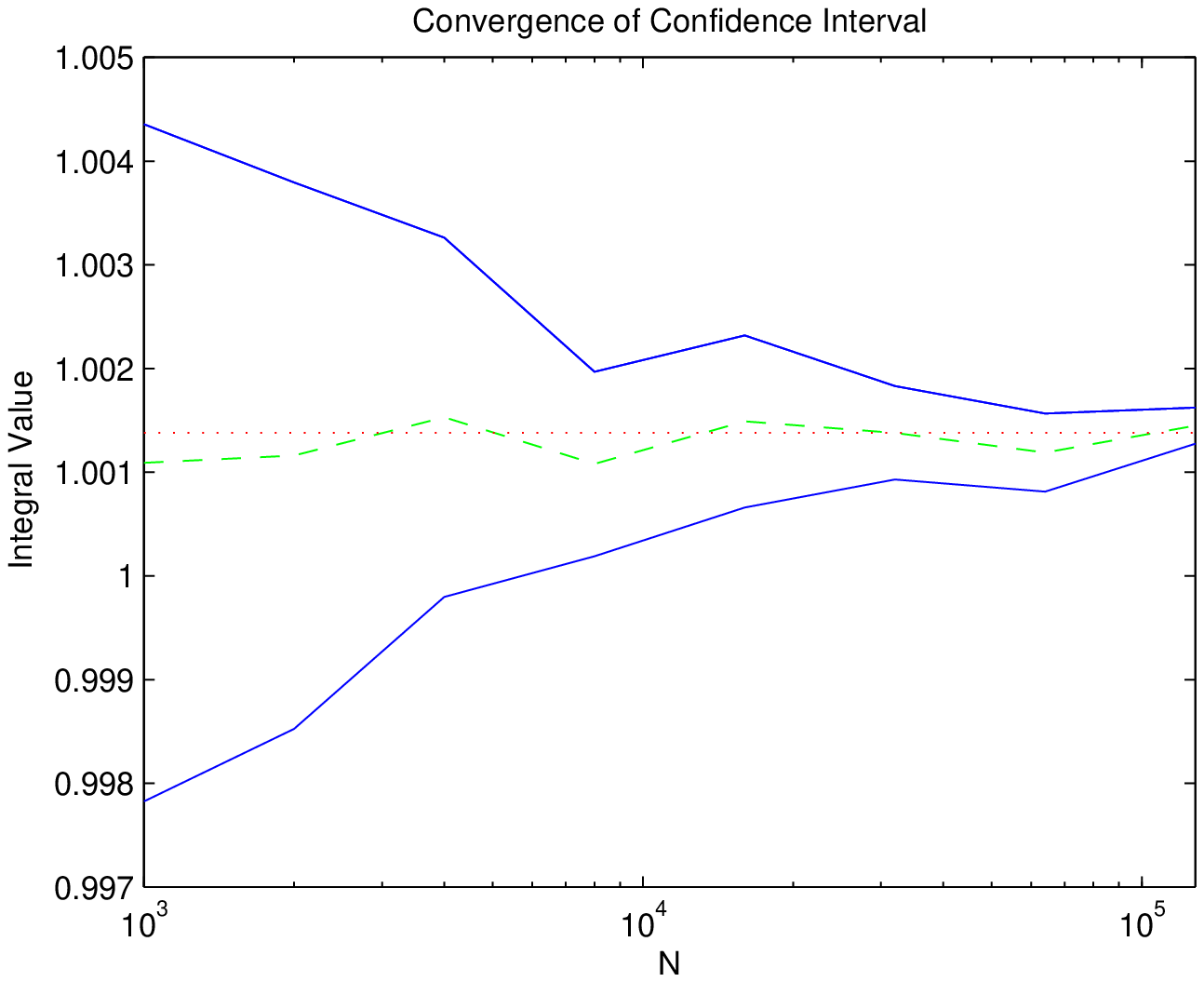}
\caption{Solid - $\mu\pm\sigma$, Dashed - $\mu$, Dotted - true solution \label{fig:confconv}}
\end{figure}

We can see from \figref{fig:diffN} that increasing $N$ has the effect of decreasing the standard deviation of the estimator.  \figref{fig:confconv} is a graphic depiction of this.  It appears in \figref{fig:diffn} and \figref{fig:diffell} that increasing $n$ and $\ell$ without changing $N$ causes no improvement in $\sigma$.  This leads us to believe that the driving force behind accuracy is $N$ for which the complexity of the algorithm increases linearly.

There are limitations to this algorithm because it requires storage of $N+1\times\ell$ terms: this is done to prevent the loss of accuracy from $N$ applications of $T$ on the elements of $S$.  Unfortunately, since there is only one point in $S$ which minimizes \eqref{optimal} there is only one pseudorandom trajectory chosen per $S$ generated.

One possible future improvement to this algorithm may be to use several trajectories whose value in \eqref{optimal} are close to optimal but not the exact minimum.  Their contribution to the estimator can be weighted according to their distance from the optimal value.  This would allow the use of multiple trajectories from the same $S$ and not require $\alpha$ versions of $S$ for $\alpha$ trajectories.

\newpage

\section*{Appendix A - Matlab Algorithms}
{\footnotesize {The first file is a script which calls the other
functions to generate pseudorandom numbers on the Julia Set for
the mapping $T(z)=z^2+1/8$.  Here is a list of important
parameters:
\begin{itemize}
\item \texttt{z0} - A repelling periodic point which is the start of the inverse iterations
\item \texttt{m} - $2^m$ Borel sets are used to cover the Julia Set
\item \texttt{zstar} - The centers of the Borel sets.  These are the 2\^{}$m$ pre-images in the set $T$\^{}$\{-m\}(z_0)$
\item \texttt{S} - The discretization of the Julia Set
\item \texttt{ell} - The number of points in $S$
\item \texttt{N} - All points in $S$ are pre-images in the set $T$\^{}$\{-(m+N)\}(z_0)$
\item \texttt{n} - Summations in the optimization equation are of length $n$
\item \texttt{alpha} - The required number of pseudorandom points
\item \texttt{h} - The Haussdorff dimension
\end{itemize}

\section*{A.1 mainscript.m}
\lstset{language=Matlab,keywordstyle=\color{blue}\bfseries,commentstyle=\color{green},basicstyle=\small,
showspaces=false,tabsize=2}
\begin{lstlisting}
% Here we consider the mapping T(z)=z^2+1/8
T=inline('x.*x+.125');
z0=fsolve(@(z) T(z)-z,.8); % approximately z0=0.85355339203135
m=8;
ell=100;
N=32000;
n=100;
alpha=30;
h=1.00735;
indmin = zeros(alpha,1);
beta_n = zeros(alpha,1);
caverage = zeros(alpha,1);
% This line makes the random number generator start with the
% same seed always.  The line below will randomize the seed.
rand('state',0);
% rand('state',sum(100*clock));

% We use inverse iteration to create the lattice.  We must
% determine the points in zstar and define the balls which
% cover the Julia Set.
zstar=z0*ones(2^m,1);
for bdec=1:2^m
  bcode=dec2bin(bdec-1,m);
  for i=1:m
    zstar(bdec)=(-1)^(bcode(i)=='1')*sqrt(zstar(bdec)-.125);
  end
end
% The function densop finds the density at zstar.
% This does not have a time limit, so if the script is
% hanging, densop is a likely source.  Both the density of
% the points in zstar and their images T(zstar) are found.
fzstar=zeros(2^m,1);
fTzstar=zeros(2^m,1);
for k=1:2:2^m
  fzstar(k)=densop(zstar(k),h);
  fzstar(k+1)=densop(zstar(k+1),h);
  fTzstar(k)=densop(T(zstar(k)),h);
  fTzstar(k+1)=fTzstar(k);
end

% Now use makelattice to form the discrete Julia Set.  The
% seeds for the inverse iteration of points on the lattice
% are randomly chosen from zstar.  The computational lattice
% is S(:,N). The final column is an extra inverse iteration
% for evaluating the fundamental equation.  It is important
% to note that the size of S is [ell,N+1] not [ell,N].

% After that we use the optimization function which will return
% the solution to the least squares problem detailed earlier.
% The function below returns imin, the index of the point in S
% which is the solution, and minival which is the residual.
% We should have beta_n/n^2<2*dn*|A_n|*LipConst as described
% in the earlier paper.
for i=1:alpha
  S=makelattice(ell,N,m,zstar);
  [imin(i),beta_n(i)]=opteval(n,S,zstar,fzstar,fTzstar,m,h);
end

% Now we test the ensemble averaging.  Any L1 function can be
% used to test the pseudorandomness of the points found by the
% optimization.  The function must be able to accept vector
% arguments, ie using .* instead of just * for multiplication.
g=inline('abs(x)');
for i=1:alpha
  caverage(i)=mean(g(S(imin(i),1:N)));
end
\end{lstlisting}

\section*{A.2 makelattice.m}
\begin{lstlisting}
function S=makelattice(ell,N,m,zstar)
% function S=makelattice(ell,N,m,zstar)
% This function makes the computational lattice that represents
% the Julia Set.  The rest of the matrix values are
% the trajectories taken backwards from a random selection of
% points on the zstar grid.  S(:,N) is the computational
% lattice, S(:,1) is ell randomly chosen points from zstar.
% S(:,N+1) is a preimage of S(:,N) which is needed for
% optimization to test if S(:,N) is in TA.
S = 2*(rand(ell,N+1)>.5)-1;
S(:,1) = zstar(ceil(2^m*rand(ell,1)));
for j=2:N+1
    S(:,j) = S(:,j).*sqrt(S(:,j-1)-.125);
end
\end{lstlisting}

\section*{A.3 densop.m}
\begin{lstlisting}
function fz=densop(z0,h)
% function fz=densop(z0,h)
% This function computes the density for the
% mapping T(z)=z^2+1/8.  It does this using a limiting sequence.
fval=[0,-1];
c=1;
cc=1;
pq(1)=z0;
pn(1)=abs(pq(1));
k=2;
while abs(fval(2)-fval(1))>1e-4
  c=[c,2^(k-1)];
  cc=cumsum(c);
  pq=[pq,zeros(1,c(k))];
  pn=[pn,zeros(1,c(k))];
  fval(1)=fval(2);
  fval(2)=0;
  for j=c(k):2:cc(k)
    pq(j)=-1^j*sqrt(pq(fix(j/2))-.125);
    pq(j+1)=-pq(j);
    pn(j)=abs(pq(j));
    pn(j+1)=pn(j);
  end
  kk=c(k);
  while kk<cc(k)
    kj=kk;
    while kj(length(kj))>2
      kj=[kj,kj(length(kj))/2-(mod(kj(length(kj)),4)>0)];
    end
    fval(2)=fval(2)+(prod(pn(kj)))^-h;
    kk=kk+2;
  end
  fval(2)=fval(2)*2^(-h*(k-1)+1);
  k=k+1;
end
fz=fval(2);
\end{lstlisting}

\section*{A.4 opteval.m}
\begin{lstlisting}
function [indmin,minival]=opteval(n,S,zstar,fzstar,fTzstar,m,h)
% function [indmin,minival]=opteval(n,S,zstar,fzstar,fTzstar,m,h)
% This needs the computational lattice, the transfer operator
% evaluated on the lattice, and the size of the Borel Set A
% around zt. n is the limit truncation which can not be greater
% than N, and h is the Hausdorff dimension.

% All this function does is run through the lattice and calculate
% the optimization equation at each point.  It returns the lowest
% value.  S is the group of trajectories which yield the
% computational lattice. Specifically S(:,N) is the lattice.
[ell,N]=size(S);
N=N-1; % Recall size(S)=[ell,N+1] although S(:,N) is the lattice.
% The radius of the balls which cover the region is related to m
delta_A=2^(-m+1);

% The lhs part of the summation will find whether either of the
% preimages of S are in A.  This is equal to asking if S is in
% T(A).  The rhs part of the summation tests whether S is in A,
% and then adds the appropriate values.  The rest is just
% evaluating the optimization equation.

Asum = zeros(ell,1);
for j=1:2^m
    spS = sparse(abs(S(:,N-n:N)-zstar(j))<delta_A);
    rhs = 2*sum((spS.*abs(S(:,N-n:N))).^h,2);
    lhs = sum((S(:,N+1-n:N+1)-zstar(j)<delta_A) + ...
              (-S(:,N+1-n:N+1)-zstar(j)<delta_A),2);
    Asum = Asum+(lhs/fTzstar(j)-rhs/fzstar(j)).^2;
end
[minival,indmin] = min(Asum/N);
\end{lstlisting}


 \bigskip
 {\bf \Large Acknowledgements}

 This work was partly supported by the NSF Grant  DMS-0620539. We
 would like to thank Xiaofan Li for helpful discussions in
 numerical error analysis. This work is partly done while J.
 Duan was a Visiting Professor at Universit\"{a}t G\"{o}ttingen and while
 M. Denker is a Distinguished Research Professor at Illinois
 Institute of Technology.


\end{document}